\documentclass{amsart}
\author{Heike Mildenberger}
\author{Dilip Raghavan}
\thanks{Second author partially supported by Grants-in-Aid for Scientific Research for JSPS
Fellows No.\ 23$\cdot$01017}
\author{Juris Stepr{\=a}ns}
\thanks{Third author partially supported by NSERC}
\date{\today}
\title{Splitting families and complete separability}
\usepackage{amssymb, amsmath, amsthm, mathrsfs, enumerate, amsfonts, latexsym, bbm}
\def\polhk#1{\setbox0=\hbox{#1}{\ooalign{\hidewidth
    \lower1.5ex\hbox{`}\hidewidth\crcr\unhbox0}}}
\newtheorem{Theorem}{Theorem}

\newtheorem{Lemma}[Theorem]{Lemma}

\newtheorem{Question}[Theorem]{Question}
\theoremstyle{definition}
\newtheorem{Def}[Theorem]{Definition}

\theoremstyle{remark}

\newcommand{\restrict}{\upharpoonright}
\newcommand{\forallbutfin}{{\forall}^{\infty}}
\newcommand{\existsinf}{{\exists}^{\infty}}
\renewcommand{\c}{\mathfrak{c}}
\renewcommand{\b}{\mathfrak{b}}
\renewcommand{\d}{{\mathfrak{d}}}
\newcommand{\s}{\mathfrak{s}}

\renewcommand{\a}{{\mathfrak{a}}}
\renewcommand{\[}{\left[}
\renewcommand{\]}{\right]}

\newcommand{\lc}{\left|}
\newcommand{\rc}{\right|}

\newcommand\ZFC{\mathrm{ZFC}}

\newcommand\PCF{\mathrm{PCF}}

\newcommand{\BS}{{\omega}^{\omega}}

\DeclareMathOperator{\otp}{otp}

\DeclareMathOperator{\dom}{dom}

\DeclareMathOperator{\cf}{cf}
\newcommand{\Pset}{\mathcal{P}}

\newcommand{\A}{{\mathscr{A}}}

\newcommand{\GG}{{\mathcal{G}}}

\newcommand{\cube}{{\[\omega\]}^{\omega}}
\newcommand{\fin}{{{\[\omega\]}^{< \omega}}}
\newcommand{\I}{{\mathcal{I}}}
\newcommand{\F}{{\mathcal{F}}}

\begin{document}
\begin{abstract}
	We answer a question from Raghavan and Stepr{\=a}ns~\cite{weaklytight} by showing that $\s = {\s}_{\omega, \omega}$. Then we use this to construct a completely separable maximal almost disjoint family under $\s \leq \a$, partially answering a question of Shelah~\cite{Sh:935}. 
\end{abstract}
\maketitle
\section{Introduction} \label{sec:intro}
	The purpose of this short note is to answer a question posed by the second and third authors in \cite{weaklytight} and to use this to solve a problem of Shelah \cite{Sh:935}. We say that two infinite subsets $a$ and $b$ of $\omega$ are \emph{almost disjoint or a.d.\@} if $a \cap b$ is finite. We say that a family $\A$ of infinite subsets of $\omega$ is \emph{almost disjoint or a.d.\@} if its members are pairwise almost disjoint. A \emph{Maximal Almost Disjoint family, or MAD family} is an infinite a.d.\ family that is not properly contained in a larger a.d.\ family.

	For an a.d.\ family $\A$, let $\I(\A)$ denote \emph{the ideal on $\omega$ generated by $\A$} -- that is, $a \in \I(\A)$ iff $\exists {a}_{0}, \dotsc, {a}_{k} \in \A \[a \; {\subset}^{\ast} \; {a}_{0} \cup \dotsb \cup {a}_{k} \]$. For any ideal $\I$ on $\omega$, ${\I}^{+}$ denotes $\Pset(\omega) \setminus \I$. An a.d.\ family $\A \subset \cube$ is said to be \emph{completely separable} if for any $b \in {\I}^{+}(\A)$, there is an $a \in \A$ with $a \subset b$. Notice that an infinite completely separable a.d.\ $\A$ must be MAD. Though the following is one of the most well-studied problems in set theory, it continues to remains open. 
\begin{Question} [Erd{\H{o}}s and Shelah~\cite{ersh}, 1972]\label{Q:compsep}
	Does there exist a completely separable MAD family $\A \subset \cube$?
\end{Question}
Progress on Question \ref{Q:compsep} was made by Balcar, Do{\v{c}}k{\'a}lkov{\'a}, and Simon who showed in a series of papers that completely separable MAD families can be constructed from any of the assumptions $\b = \d$, $\s = {\omega}_{1}$, or $\d \leq \a$. See \cite{BS}, \cite{disjoint}, and \cite{simon1} for this work. Then Shelah~\cite{Sh:935} recently showed that the existence of completely separable MAD families is \emph{almost} a theorem of $\ZFC$. His construction is divided into three cases. The first case is when $\s < \a$ and he shows on the basis of $\ZFC$ alone that a completely separable MAD family can be constructed in this case. The second and third cases are when $\s = \a$ and $\a < \s$ respectively and Shelah shows that a completely separable MAD family can be constructed in these cases \emph{provided} that certain $\PCF$ type hypotheses are satisfied. More precisely, he shows that there is a completely separable MAD family when $\s = \a$ and $U(\s)$ holds, or when $\a < \s$ and $P(\s, \a)$ holds.
\begin{Def} \label{def:UandP}
	For a cardinal $\kappa > \omega$, $U(\kappa)$ is the following principle. There is a sequence $\langle {u}_{\alpha}: \omega \leq \alpha < \kappa \rangle$ such that
	\begin{enumerate}
		\item
			${u}_{\alpha} \subset \alpha$ and $\lc {u}_{\alpha} \rc = \omega$	
		\item
			$\forall X \in {\[\kappa\]}^{\kappa} \exists \omega \leq \alpha < \kappa\[\lc {u}_{\alpha} \cap X \rc = \omega\]$.
	\end{enumerate}
	For cardinals $\kappa > \lambda > \omega$, $P(\kappa, \lambda)$ says that there is a sequence $\langle {u}_{\alpha}: \omega \leq \alpha < \kappa \rangle$ such that
	\begin{enumerate}
		\item[(3)]
			${u}_{\alpha} \subset \alpha$ and $\lc {u}_{\alpha} \rc = \omega$	
		\item[(4)]
			for each $X \subset \kappa$, if $X$ is bounded in $\kappa$ and $\otp{(X)} = \lambda$, then $\exists \omega \leq \alpha < \sup{(X)}\[\lc {u}_{\alpha} \cap X \rc = \omega\]$.
	\end{enumerate}
\end{Def}
It is easy to see that both $U(\s)$ and $P(\s, \a)$ are satisfied when $\s < {\aleph}_{\omega}$, so in particular, the existence of a completely separable MAD family is a theorem of $\ZFC$ when $\c < {\aleph}_{\omega}$. Shelah~\cite{Sh:935} asked whether all uses of $\PCF$ type hypotheses can be eliminated from the second and third cases.

	The second and third authors modified the techniques of Shelah~\cite{Sh:935} in order to treat MAD families  with few partitioners in \cite{weaklytight} (see the introduction there). In that paper they introduced a cardinal invariant ${\s}_{\omega, \omega}$, which is a variation of the splitting number $\s$. They showed that if ${\s}_{\omega, \omega} \leq \b$, then there is a weakly tight family. Recall that an a.d.\ family $\A \subset \cube$ is called \emph{weakly tight} if for every countable collection $\{{b}_{n}: n \in \omega\} \subset {\I}^{+}(\A)$, there is $a \in \A$ such that $\existsinf n \in \omega \[\lc {b}_{n} \cap a \rc = \omega \]$. The question of whether $\s = {\s}_{\omega, \omega}$ was raised in \cite{weaklytight}, and the authors pointed out that an affirmative answer to this question could help eliminate the use of $\PCF$ type hypotheses from the second case of Shelah's construction.

	In this paper we answer this question from \cite{weaklytight} by proving that $\s = {\s}_{\omega, \omega}$. We then use this information to partially answer the question from Shelah~\cite{Sh:935}. We show that the second case can be done without any additional hypothesis. So it is a theorem of $\ZFC$ alone that a completely separable MAD family exists when $\s \leq \a$. We give a single construction from this assumption, so Shelah's first and second cases are unified into a single case.

	The question of whether the hypothesis $P(\s, \a)$ can be eliminated from the case when $\a < \s$ remains open.       
\section{$\s = {\s}_{\omega, \omega}$} \label{sec:sequalssomegaomega}
	In this section we answer Question 21 from \cite{weaklytight} by showing that $\s = {\s}_{\omega, \omega}$. For a set $x \subset \omega$, ${x}^{0}$ is used to denote $x$ and ${x}^{1}$ is used to denote $\omega \setminus x$. This notation will be used in the next section also. Recall the following definitions.
\begin{Def} \label{def:sandsomegaomega}
 	For $x, a \in \Pset(\omega)$, $x$ \emph{splits} $a$ if $\lc {x}^{0} \cap a \rc = \lc {x}^{1} \cap a \rc = \omega$. $\F \subset \Pset(\omega)$ is called a \emph{splitting family} if $\forall a \in \cube \exists x \in \F\[x \ \text{splits} \ a\]$. $\F \subset \Pset(\omega)$ is said to be \emph{$(\omega, \omega)$-splitting} if for each countable collection $\{{a}_{n}: n \in \omega\} \subset \cube$, there exists $ x \in \F$ such that $\existsinf n \in \omega \[\lc {x}^{0} \cap {a}_{n} \rc = \omega \]$ and $\existsinf n \in \omega \[\lc {x}^{1} \cap {a}_{n} \rc = \omega \]$. Define	
	\begin{align*}
		&\s = \min\{\lc \F \rc: \F \subset \Pset(\omega) \wedge \F \ \text{is a splitting family}\}\\
		&{\s}_{\omega, \omega} = \min\{\lc \F \rc: \F \subset \Pset(\omega) \wedge \F \ \text{is} \ (\omega, \omega)-\text{splitting}\}.
	\end{align*}
\end{Def}
Obviously every $(\omega, \omega)$-splitting family is a splitting family. So $\s \leq {\s}_{\omega, \omega}$. It was shown in Theorem 13 of \cite{weaklytight} that if $\s < \b$, then $\s = {\s}_{\omega, \omega}$. We reproduce that result here for the reader's convenience.
\begin{Lemma} [Theorem 13 of \cite{weaklytight}]\label{lem:whenbisbig}
	If $\s < \b$, then $\s = {\s}_{\omega, \omega}$.
\end{Lemma}
\begin{proof}
	Let $\langle {e}_{\alpha}: \alpha < \kappa \rangle$ witness that $\kappa = \s$. Suppose $\{{b}_{n}: n \in \omega \} \subset \cube$ is a countable collection such that $\forall \alpha < \kappa \exists i \in 2 \forallbutfin n \in \omega \[{b}_{n} \; {\subset}^{\ast} \; {e}^{i}_{\alpha}\]$. By shrinking them if necessary we may assume that ${b}_{n} \cap {b}_{m} = 0$ whenever $n \neq m$. Now, for each $\alpha < \kappa$ define ${f}_{\alpha} \in \BS$ as follows. We know that there is a unique ${i}_{\alpha} \in 2$ such that there is a ${k}_{\alpha} \in \omega$ such that $\forall n \geq {k}_{\alpha} \[\lc {b}_{n} \cap {e}^{{i}_{\alpha}}_{\alpha} \rc < \omega \]$. We define ${f}_{\alpha}(n) = \max\left( {b}_{n} \cap {e}^{{i}_{\alpha}}_{\alpha} \right)$ if $n \geq {k}_{\alpha}$, and ${f}_{\alpha}(n) = 0$ if $n < {k}_{\alpha}$. As $\kappa < \b$, there is a $f \in \BS$ with $f \; {}^{\ast}{>} \; {f}_{\alpha}$ for each $\alpha < \kappa$. Now, for each $n \in \omega$, choose ${l}_{n} \in {b}_{n}$ with ${l}_{n} \geq f(n)$. Since the ${b}_{n}$ are pairwise disjoint, $c = \{{l}_{n}: n \in \omega\} \in \cube$. So by definition of $\s$, there is $\alpha < \kappa$ such that $\lc c \cap {e}^{0}_{\alpha} \rc = \lc c \cap {e}^{1}_{\alpha} \rc = \omega$. In particular, $c \cap {e}^{{i}_{\alpha}}_{\alpha}$ is infinite. But we know that there is an ${m}_{\alpha} \in \omega$ such that $\forall n \geq {m}_{\alpha}\[{f}_{\alpha}(n) < f(n)\]$. So there exists $n \geq \max\{{m}_{\alpha}, {k}_{\alpha}\}$ with ${l}_{n} \in {b}_{n} \cap {e}^{{i}_{\alpha}}_{\alpha}$. But this is a contradiction because ${l}_{n} \leq {f}_{\alpha}(n) < f(n)$.
\end{proof}
In the case when $\b \leq \s$ it turns out that $\s = {\s}_{\omega, \omega}$ can still be proved by considering the following notion appearing in \cite{KW}.  
\begin{Def} \label{def:block}
$\F$ is called \emph{block-splitting} if given any partition $\langle
a_n : n \in \omega \rangle$ of $\omega$ into finite sets there is a set 
$x \in \F$ such that there are infinitely many $n$ with $a_n \subset x$ 
and there are infinitely many $n$ with $a_n \cap x = 0$.
\end{Def}
It was proved by Kamburelis and W{\polhk{e}}glorz~\cite{KW} that the least size of a block splitting family is $\max\{\b, \s\}$. Therefore, when $\b \leq \s$, there is a block splitting family of size $\s$.
\begin{Theorem} \label{thm:sequalssomegaomega}
	$\s = {\s}_{\omega, \omega}$.
\end{Theorem}
\begin{proof}
	In view of Lemma \ref{lem:whenbisbig}, we may assume that $\b \leq \s$. By results of Kamburelis and W{\polhk{e}}glorz~\cite{KW} fix $\langle {x}_{\alpha}: \alpha < \s \rangle \subset \Pset(\omega)$, a block splitting family. We show that $\langle {x}_{\alpha}: \alpha < \s \rangle$ is an $\left( \omega, \omega \right)$-splitting family. Let $\{{a}_{n}: n \in \omega\} \subset \cube$ be given. For $n \in \omega$, define ${s}_{n} \in \fin$ as follows. Suppose $\langle {s}_{i}: i < n \rangle$ have been defined. Put $s = {\bigcup}_{i < n}{{s}_{i}}$. Put ${s}_{n} = \{\min(\omega \setminus s)\} \cup \{\min({a}_{i} \setminus s): i \leq n\}$. Note that $\langle {s}_{n}: n \in \omega \rangle$ is a partition of $\omega$ into finite sets and that $\forall i \in \omega \forallbutfin n \in \omega \[{s}_{n} \cap {a}_{i} \neq 0\]$. Now choose $\alpha < \s$ such that $\existsinf n \in \omega \[{s}_{n} \subset {x}^{0}_{\alpha}\]$ and $\existsinf n \in \omega \[{s}_{n} \subset {x}^{1}_{\alpha}\]$. So for each $i \in \omega$, $\existsinf n \in \omega \[{s}_{n} \cap {a}_{i} \cap {x}^{0}_{\alpha} \neq 0 \]$ and $\existsinf n \in \omega \[{s}_{n} \cap {a}_{i} \cap {x}^{1}_{\alpha} \neq 0 \]$. Since the ${s}_{n}$ are pairwise disjoint, it follows that $\lc {a}_{i} \cap {x}^{0}_{\alpha} \rc = \lc {a}_{i} \cap {x}^{1}_{\alpha}\rc = \omega$, for each $i \in \omega$. 
\end{proof}
\section{Constructing a completely separable MAD family from $\s \leq \a$} \label{sec:main}
As $\s = {\s}_{\omega, \omega}$ and as every $\left( \omega, \omega \right)$-splitting family is also a splitting family, fix once and for all a sequence $\langle {x}_{\alpha}: \alpha < \kappa \rangle$ witnessing that $\kappa = \s = {\s}_{\omega, \omega}$. We will construct a completely separable MAD family assuming that $\kappa \leq \a$. The construction closely follows the proof of Lemma 8 in \cite{weaklytight}, which in turn is based on Shelah \cite{Sh:935}. An important point of the construction is that if $\A$ is an arbitrary a.d.\ family and $b \in {\I}^{+}(\A)$, then every $\left( \omega, \omega\right)$-splitting family contains an element which splits $b$ into two \emph{positive} pieces.  
\begin{Lemma} \label{lem:positivesplitting}
Let $\A \subset \cube$ be any a.d.\ family.
Suppose $b \in {\I}^{+}(\A)$.
Then there is $\alpha < \kappa$ such that $b \cap {x}^{0}_{\alpha} \in {\I}^{+}(\A)$ and $b \cap {x}^{1}_{\alpha} \in {\I}^{+}(\A)$. 
\end{Lemma}
\begin{proof}
	See proof of Lemma 7 of \cite{weaklytight}.
\end{proof}
At a stage $\delta < \c$, an a.d.\ family ${\A}_{\delta} = \langle {a}_{\alpha}: \alpha < \delta \rangle \subset \cube$ is given. 
Moreover we assume that there is also a family $\langle {\sigma}_{\alpha}: \alpha < \delta \rangle \subset {2}^{< \kappa}$ such that for each $\alpha < \delta$, $\forall \xi < \dom({\sigma}_{\alpha})\[{a}_{\alpha} \; {\subset}^{\ast} \; {x}^{{\sigma}_{\alpha}(\xi)}_{\xi}\]$.
We say that ${\sigma}_{\alpha}$ is \emph{the node associated with} ${a}_{\alpha}$.
The next lemma says that under the assumption $\kappa \leq \a$, such an a.d.\ family must be ``nowhere maximal'', which is of course a property that we need to maintain in order to end up with a completely separable MAD family.
\begin{Def} \label{def:ieta}
Let $\eta \in {2}^{< \kappa}$.
Define ${\I}_{\eta} = \left\{a \in \Pset(\omega): \forall \xi < \dom(\eta)\[a \; {\subset}^{\ast} \; {x}^{\eta(\xi)}_{\xi} \]\right\}$.
\end{Def}
\begin{Lemma} [Main Lemma]\label{lem:main}
Assume $\kappa \leq \a$.
Let $\delta < \c$.
Suppose that ${\A}_{\delta} = \langle {a}_{\alpha}: \alpha < \delta\rangle$ and $\langle {\sigma}_{\alpha}: \alpha < \delta \rangle$ are as above.
Assume also that $\forall \alpha, \beta < \delta \[\alpha \neq \beta \implies {\sigma}_{\alpha} \neq {\sigma}_{\beta}\]$.
Let $b \in {\I}^{+}({\A}_{\delta})$.
Then there exist $a \in {\[b\]}^{\omega}$ and $\sigma \in {2}^{< \kappa}$ such that 
	\begin{enumerate}
		\item
			$\forall \alpha < \delta \[\lc a \cap {a}_{\alpha} \rc < \omega\]$.
		\item
			for each $\alpha < \delta$, $\sigma \not\subset {\sigma}_{\alpha}$ and $a \in {I}_{\sigma}$.
	\end{enumerate}
\end{Lemma}
\begin{proof}
Applying Lemma \ref{lem:positivesplitting}, let ${\alpha}_{0} < \kappa$ be least such that $b \cap {x}^{0}_{{\alpha}_{0}} \in {\I}^{+}({\A}_{\delta})$ and $b \cap {x}^{1}_{{\alpha}_{0}} \in {\I}^{+}({\A}_{\delta})$.
Define ${\tau}_{0} \in {2}^{{\alpha}_{0}}$ by stipulating that 
	\begin{align*}
		\forall \xi < {\alpha}_{0} \forall i \in 2 \[{\tau}_{0}(\xi) = i \leftrightarrow b \cap {x}^{i}_{\xi} \in {\I}^{+}({\A}_{\delta})\].
	\end{align*} 
By choice of ${\alpha}_{0}$ and by the hypothesis that $b \in {\I}^{+}({\A}_{\delta})$, ${\tau}_{0}$ is well defined.
Now, construct two sequences $\langle {\alpha}_{s}: s \in {2}^{< \omega}\rangle \subset \kappa$ and $\langle {\tau}_{s}: s \in {2}^{< \omega}\rangle \subset {2}^{< \kappa}$ such that
	\begin{enumerate}
		\item[(3)]
			$\forall s \in {2}^{< \omega} \forall i \in 2 \[{\alpha}_{s} = \dom({\tau}_{s}) \wedge {\alpha}_{{s}^{\frown}{\langle i \rangle}} > {\alpha}_{s} \wedge {\tau}_{{s}^{\frown}{\langle i \rangle}} \supset {{\tau}_{s}}^{\frown}{\langle i \rangle}\]$.
		\item[(4)]
			for each $s \in {2}^{< \omega}$ and for each $\xi < {\alpha}_{s}$, ${x}^{1 - {\tau}_{s}(\xi)}_{\xi} \cap b \cap \left( {\bigcap}_{t \subsetneq s}{{x}^{{\tau}_{s}({\alpha}_{t})}_{{\alpha}_{t}}}\right) \notin {\I}^{+}({\A}_{\delta})$.
Here when $s = 0$, ${\bigcap}_{t \subsetneq s}{{x}^{{\tau}_{s}({\alpha}_{t})}_{{\alpha}_{t}}}$ is taken to be $\omega$. 
		\item[(5)]
			for each $s \in {2}^{< \omega}$, both ${x}^{0}_{{\alpha}_{s}} \cap b \cap \left( {\bigcap}_{t \subsetneq s}{{x}^{{\tau}_{s}({\alpha}_{t})}_{{\alpha}_{t}}} \right) \in {\I}^{+}({\A}_{\delta})$ and ${x}^{1}_{{\alpha}_{s}} \cap b \cap \left( {\bigcap}_{t \subsetneq s}{{x}^{{\tau}_{s}({\alpha}_{t})}_{{\alpha}_{t}}} \right) \in {\I}^{+}({\A}_{\delta})$.
	\end{enumerate}
${\alpha}_{0}$ and ${\tau}_{0}$ are already defined.
Suppose that ${\alpha}_{s}$ and ${\tau}_{s}$ are given.
By (5) for each $i \in 2$, ${x}^{i}_{{\alpha}_{s}} \cap b \cap \left({\bigcap}_{t \subsetneq s}{{x}^{{\tau}_{s}({\alpha}_{t})}_{{\alpha}_{t}}}\right) \in {\I}^{+}({\A}_{\delta})$.
 Apply Lemma \ref{lem:positivesplitting} to let ${\alpha}_{{s}^{\frown}{\langle i \rangle}}$ be the least $\alpha < \kappa$ such that both ${x}^{i}_{{\alpha}_{s}} \cap b \cap \left({\bigcap}_{t \subsetneq s}{{x}^{{\tau}_{s}({\alpha}_{t})}_{{\alpha}_{t}}}\right) \cap {x}^{0}_{\alpha}$ and ${x}^{i}_{{\alpha}_{s}} \cap b \cap \left({\bigcap}_{t \subsetneq s}{{x}^{{\tau}_{s}({\alpha}_{t})}_{{\alpha}_{t}}}\right) \cap {x}^{1}_{\alpha}$ are in ${\I}^{+}({\A}_{\delta})$.
Again define ${\tau}_{{s}^{\frown}{\langle i \rangle}} \in {2}^{{\alpha}_{{s}^{\frown}{\langle i \rangle}}}$ by stipulating that
	\begin{align*}
		\forall \xi < {\alpha}_{{s}^{\frown}{\langle i \rangle}} \forall j \in 2 \[{\tau}_{{s}^{\frown}{\langle i \rangle}}(\xi) = j \leftrightarrow {x}^{i}_{{\alpha}_{s}} \cap b \cap \left({\bigcap}_{t \subsetneq s}{{x}^{{\tau}_{s}({\alpha}_{t})}_{{\alpha}_{t}}}\right) \cap {x}^{j}_{\xi} \in {\I}^{+}({\A}_{\delta})\]
	\end{align*}
${\tau}_{{s}^{\frown}{\langle i \rangle}}$ is well defined because ${x}^{i}_{{\alpha}_{s}} \cap b \cap \left({\bigcap}_{t \subsetneq s}{{x}^{{\tau}_{s}({\alpha}_{t})}_{{\alpha}_{t}}}\right) \in {\I}^{+}({\A}_{\delta})$ and because of the choice of ${\alpha}_{{s}^{\frown}{\langle i \rangle}}$.
Now, for each $\xi < {\alpha}_{s}$, ${x}^{i}_{{\alpha}_{s}} \cap b \cap \left({\bigcap}_{t \subsetneq s}{{x}^{{\tau}_{s}({\alpha}_{t})}_{{\alpha}_{t}}}\right) \subset b \cap \left({\bigcap}_{t \subsetneq s}{{x}^{{\tau}_{s}({\alpha}_{t})}_{{\alpha}_{t}}}\right)$ and, by (4), $b \cap \left({\bigcap}_{t \subsetneq s}{{x}^{{\tau}_{s}({\alpha}_{t})}_{{\alpha}_{t}}}\right) \cap {x}^{1 - {\tau}_{s}(\xi)}_{\xi} \notin {\I}^{+}({\A}_{\delta})$.
It follows that ${\alpha}_{{s}^{\frown}{\langle i \rangle}} \geq {\alpha}_{s}$ and that for each $\xi < {\alpha}_{s}$, ${\tau}_{s}(\xi) = {\tau}_{{s}^{\frown}{\langle i \rangle}}(\xi)$.
Next, since ${x}^{i}_{{\alpha}_{s}} \cap b \cap \left({\bigcap}_{t \subsetneq s}{{x}^{{\tau}_{s}({\alpha}_{t})}_{{\alpha}_{t}}}\right) \cap {x}^{1 - i}_{{\alpha}_{s}} = 0$, ${\alpha}_{{s}^{\frown}{\langle i \rangle}} > {\alpha}_{s}$, and ${\tau}_{{s}^{\frown}{\langle i \rangle}} \supset {{\tau}_{s}}^{\frown}{\langle i \rangle}$.
Now, it is clear that (4) and (5) hold for ${s}^{\frown}{\langle i \rangle}$.
This completes the construction of $\langle {\alpha}_{s}: s \in {2}^{< \omega}\rangle$ and $\langle {\tau}_{s}: s \in {2}^{< \omega}\rangle$.

	For each $f \in {2}^{\omega}$, put ${\alpha}_{f} = \sup\left\{{\alpha}_{f \restrict n}: n \in \omega\right\}$ and ${\tau}_{f} = {\bigcup}_{n \in \omega}{{\tau}_{f \restrict n}}$.
As $\kappa = \s$, $\cf(\kappa) > \omega$.
Therefore, ${\alpha}_{f} < \kappa$.
Note that ${\tau}_{f} \in {2}^{{\alpha}_{f}}$.
Also, if $f, g \in {2}^{\omega}$, $f \neq g$, and $n \in \omega$ is least such that $f(n) \neq g(n)$, then ${\tau}_{f} \supset {{\tau}_{s}}^{\frown}{\langle i \rangle}$ and ${\tau}_{g} \supset {{\tau}_{s}}^{\frown}{\langle 1 - i\rangle}$, where $s = f \restrict n = g \restrict n$ and $i \in 2$.
So there cannot be $\alpha < \delta$ such that both ${\tau}_{f} \subset {\sigma}_{\alpha}$ and ${\tau}_{g} \subset {\sigma}_{\alpha}$ hold.
Therefore, it is possible to find $f \in {2}^{\omega}$ such that ${\tau}_{f} \not\in \{\sigma \in {2}^{< \kappa}: \exists \alpha < \delta \[\sigma \subset {\sigma}_{\alpha} \] \}$. Fix such $f$ and for each $n \in \omega$, define ${e}_{n}$ to be $b \cap \left({\bigcap}_{m < n}{{x}^{{\tau}_{f}\left({\alpha}_{f \restrict m}\right)}_{{\alpha}_{f \restrict m}}} \right)$. By (5) each ${e}_{n} \in {\I}^{+}({\A}_{\delta})$. Moreover, ${e}_{n + 1} \subset {e}_{n} \subset b$. Therefore, by a standard argument, there is $e \in {\[b\]}^{\omega} \cap {\I}^{+}({\A}_{\delta})$ such that $\forall n \in \omega \[e \; {\subset}^{\ast} \; {e}_{n}\]$.

	Now, suppose $\xi < {\alpha}_{f}$. Since for all $n \in \omega$, ${\alpha}_{f \restrict n + 1} > {\alpha}_{f \restrict n}$, it follows that $\xi < {\alpha}_{f \restrict n}$ for some  $n$. By (4) applied to $s = f \restrict n$, ${x}^{1 - {\tau}_{f}(\xi)}_{\xi} \cap {e}_{n} \notin {\I}^{+}({\A}_{\delta})$. Since $e \; {\subset}^{\ast} \; {e}_{n}$, ${x}^{1 - {\tau}_{f}(\xi)}_{\xi} \cap e \notin {\I}^{+}({\A}_{\delta})$. Thus we conclude that $\forall \xi < {\alpha}_{f}\[{x}^{1 - {\tau}_{f}(\xi)}_{\xi} \cap e \notin {\I}^{+}({\A}_{\delta})\]$. So for each $\xi < {\alpha}_{f}$, fix ${F}_{\xi} \in {\[\delta\]}^{< \omega}$ such that
	\begin{align*}
		\left({x}^{1 - {\tau}_{f}(\xi)}_{\xi} \cap e \right) \; {\subset}^{\ast} \; \left( {\bigcup}_{\alpha \in {F}_{\xi}}{{a}_{\alpha}} \right)
	\end{align*}  
Now, put $\F = {\bigcup}_{\xi < {\alpha}_{f}}{{F}_{\xi}}$ and $\GG = \{\alpha < \delta: {\sigma}_{\alpha} \subset {\tau}_{f}\}$. Note that $\lc \F  \cup \GG \rc < \kappa \leq \a$ because of the assumption that $\forall \alpha, \beta < \delta \[\alpha \neq \beta \implies {\sigma}_{\alpha} \neq {\sigma}_{\beta} \]$. Since $e \in {\I}^{+}({\A}_{\delta})$, there is $a \in {\[e\]}^{\omega}$ such that $\forall \alpha \in \F \cup \GG \[\lc a \cap {a}_{\alpha} \rc < \omega\]$. Note that for each $\xi < {\alpha}_{f}$, ${x}^{1 - {\tau}_{f}(\xi)}_{\xi} \cap a$ is finite. Thus putting $\sigma = {\tau}_{f}$, we have that $\forall \alpha < \delta \[\sigma \not\subset  {\sigma}_{\alpha}\]$ and $a \in {I}_{\sigma}$. In order to finish the proof, it is enough to check that $\forall \alpha < \delta \[\lc {a}_{\alpha} \cap a \rc < \omega\]$. 

	Fix $\alpha < \delta$. If $\alpha \in \GG$, then $\lc a \cap {a}_{\alpha} \rc < \omega$ simply by choice of $a$. Suppose $\alpha \notin \GG$. Then there must be $\xi \in \dom({\sigma}_{\alpha}) \cap {\alpha}_{f}$ such that ${\sigma}_{\alpha}(\xi) = 1 - {\tau}_{f}(\xi)$. But since ${a}_{\alpha} \; {\subset}^{\ast} \; {x}^{{\sigma}_{\alpha}(\xi)}_{\xi}$ and $a \cap {x}^{1 - {\tau}_{f}(\xi)}_{\xi}$ is finite, it follows that $\lc a \cap {a}_{\alpha} \rc < \omega$.
\end{proof}
\begin{Theorem} \label{thm:main}
	If $\s \leq \a$, then there is a completely separable MAD family.
\end{Theorem}
\begin{proof}
	Fix an enumeration $\langle {b}_{\alpha}: \alpha < \c \rangle$ of $\cube$. Let $\langle {x}_{\alpha}: \alpha < \kappa \rangle$ witness $\kappa = \s = {\s}_{\omega, \omega}$. Build two sequences $\langle {a}_{\delta}: \delta < \c \rangle$ and $\langle {\sigma}_{\delta}: \delta < \c \rangle$ such that the following hold.
	\begin{enumerate}
		\item
			for each $\delta < \c$, ${a}_{\delta} \in \cube$, ${\sigma}_{\delta} \in {2}^{< \kappa}$, and ${a}_{\delta} \in {I}_{{\sigma}_{\delta}}$
		\item
			$\forall \gamma, \delta < \c \[\gamma \neq \delta \implies \left( \lc {a}_{\gamma} \cap {a}_{\delta} \rc < \omega \wedge {\sigma}_{\gamma} \neq {\sigma}_{\delta} \right) \]$
		\item
			for each $\delta < \c$, if ${b}_{\delta} \in {\I}^{+}({\A}_{\delta})$, then ${a}_{\delta} \subset {b}_{\delta}$, where ${\A}_{\delta} = \{{a}_{\alpha}: \alpha < \delta \}$
	\end{enumerate} 
	Note that if we succeed in this, then ${\A}_{\c} = \{{a}_{\delta}: \delta < \c \}$ will be completely separable for given any $b \in {\I}^{+}({\A}_{\c})$, $b$ is in ${\I}^{+}({\A}_{\delta})$ for every $\delta < \c$ and so there is a $\delta < \c$ where ${b}_{\delta} = b$ and ${b}_{\delta} \in {\I}^{+}({\A}_{\delta})$, whence by (3), ${a}_{\delta} \subset b$.
	
	At a stage $\delta < \c$ suppose $\langle {a}_{\alpha}: \alpha < \delta \rangle$ and $\langle {\sigma}_{\alpha}: \alpha < \delta \rangle$ are given. If ${b}_{\delta} \in {\I}^{+}({\A}_{\delta})$, then let $b = {b}_{\delta}$, else let $b = \omega$. In either case, the hypotheses of Lemma \ref{lem:main} are satisfied. So find ${a}_{\delta} \in {\[b\]}^{\omega}$ and ${\sigma}_{\delta} \in {2}^{< \kappa}$ such that
	\begin{enumerate}
		\item[(4)]
			$\forall \alpha < \delta \[\lc {a}_{\delta} \cap {a}_{\alpha } \rc < \omega \]$
		\item[(5)]
			for each $\alpha < \delta$, ${\sigma}_{\delta} \not\subset {\sigma}_{\alpha}$ and ${a}_{\delta} \in {I}_{{\sigma}_{\delta}}$.
	\end{enumerate}
It is clear that ${a}_{\delta}$ and ${\sigma}_{\delta}$ are as needed.
\end{proof}
\bibliographystyle{amsplain}
\bibliography{Bibliography}
\end{document}